\documentstyle[12pt]{article}
\newtheorem{theorem}{Theorem}\newtheorem{lemma}[theorem]{Lemma}

\def\pont{\hspace{-6pt}{\bf.\ \,}} 
\def\bs{\bigskip}
\newcommand{\text}[1]{\quad\mbox{#1}\quad}
\def\beq{\begin{equation}}\def\eeq{\end{equation}}
\def\beqn{\begin{eqnarray}}\def\eeqn{\end{eqnarray}}
\def\beqno{\begin{eqnarray*}}\def\eeqno{\end{eqnarray*}}
\def\ds{\displaystyle}\def\eps{\varepsilon}\def\ffi{\varphi}
\def\qed{\ifhmode\unskip\nobreak\fi\quad\ifmmode\Box\else$\Box$\fi}
\def\({\mbox{$($}}\def\){\mbox{$)$}}
\title{A quantitative version of the Blow-up Lemma}
\author{G\'{a}bor N. S\'ark\"ozy\thanks{Research supported in part by
OTKA Grant No. K104343.}\\
\small Alfr\'ed R\'enyi Institute of Mathematics\\[-0.8ex]
\small Hungarian Academy of Sciences\\[-0.8ex]
\small Budapest, P.O. Box 127\\[-0.8ex]
\small Budapest, Hungary, H-1364\\
and\\
\small Computer Science Department\\[-0.8ex]
\small Worcester Polytechnic Institute\\[-0.8ex]
\small Worcester, MA, USA 01609\\[-0.8ex]
\small \texttt{gsarkozy@cs.wpi.edu}}

\date{}
\begin{document}
\parindent0pt
\maketitle

\begin{abstract}
In this paper we give a quantitative version of the Blow-up Lemma.
\end{abstract}

\section{Introduction}
\subsection{Notations and definitions}

All graphs are simple, that is, they have no loops or multiple
edges. $v(G)$ is the number of vertices in $G$ (order), $e(G)$ is
the number of edges in $G$ (size). $deg(v)$ (or $deg_G(v)$) is the
degree of vertex $v$ (within the graph $G$), and $deg(v,Y)$ (or
$deg_G(v,Y)$) is the number of neighbors of $v$ in $Y$. $\delta(G)$
and $\Delta(G)$ are the minimum degree and the maximum degree of
$G$. $N(x)$ (or $N_G(x)$) is the set of neighbors of the vertex
$x$, and $e(X,Y)$ is the number of edges between $X$ and $Y$. A
bipartite graph $G$ with color-classes $A$ and $B$ and edges $E$
will sometimes be written as $G=(A,B,E)$. For disjoint $X,Y$, we
define the {\bf density}
$$d(X,Y)={e(X,Y)\over|X|\cdot|Y|}\,.$$
The density of a bipartite graph $G=(A,B,E)$ is the number
$$d(G)=d(A,B)={|E|\over|A|\cdot|B|}\,.$$
For two disjoint subsets $A,B$ of $V(G)$,
the bipartite graph with vertex set $A\cup B$ which has
all the edges of $G$ with one endpoint in $A$
and the other in $B$ is called the pair $(A,B)$.

A pair $(A,B)$ is $\eps$-{\bf regular}
if for every $X\subset A$ and $Y\subset B$ satisfying
$$|X|>\eps|A|\text{and}|Y|>\eps|B|$$
we have $$|d(X,Y)-d(A,B)|<\eps.$$ A pair $(A,B)$ is
$(\eps,d,\delta)${\bf-super-regular} if it is $\eps$-regular with
density at least $d$ and furthermore,$$deg(a)\geq\delta|B|\text{for
all}a\in A,$$
$$\text{and}deg(b)\geq\delta|A|\text{for all}b\in B.$$

$H$ is {\bf embeddable} into $G$ if $G$ has a subgraph isomorphic to $H$,
that is, if there is a one-to-one map (injection)
$\ffi:V(H)\to V(G)$ such that $\{x,y\}\in E(H)$ implies
$\{\ffi(x),\ffi(y)\}\in E(G)$.

\subsection{A quantitative version of the Blow-up Lemma}

The Blow-up Lemma \cite{KSSz2, KSSz3}  has been a successful tool in
extremal graph theory. There are now at least four new proofs for
the Blow-up Lemma since the original appeared; an algorithmic proof
\cite{KSSz3}, a hypergraph-packing approach \cite{RRT}, a proof
based on counting perfect matchings in (Szemer\'{e}di-) regular
graphs \cite{RR}, and its constructive version in \cite{RRW}. Very
recently the Blow-up Lemma has been generalized to hypergraphs by
Keevash \cite{KEE} and to $d$-arrangeable graphs by B\"ottcher,
Kohayakawa, Taraz, and W\"urfl \cite{BKTW}. The Blow-up Lemma has
been applied in numerous papers (see e.g. \cite{ARR, GRR, KOR, KS,
KRSSz, KO1, KO2, LRS, RR, RRT, RRW}). See also the discussion on the
Regularity Lemma and the Blow-up Lemma on pages 803-804 in the
Handbook of Graph Theory \cite{B2} or the survey paper \cite{KO3}.

In either of our proofs \cite{KSSz2, KSSz3}, the dependence of the
parameters was not computed explicitly. In this paper we give a
quantitative version, i.e. we compute explicitly the parameters.

\begin{theorem}{\bf (A quantitative version of the Blow-up Lemma)}
\label{blowup}\pont There exists an absolute constant $C$ such that,
given a graph $R$ of order $r\geq 2$ and positive parameters $d$,
$\delta$, and $\Delta$, for any $0< \eps < \left( \frac{\delta
d^{\Delta}}{r\Delta}\right)^C$ the following holds. Let $N$ be an
arbitrary positive integer, and let us replace the vertices of $R$
with pairwise disjoint $N$-sets $V_1,V_2,\ldots,V_r$ \(blowing up\).
We construct two graphs on the same vertex-set $V=\cup V_i$. The
graph $R(N)$ is obtained by replacing all edges of $R$ with copies
of the complete bipartite graph $K_{N,N}$, and a sparser graph $G$
is constructed by replacing the edges of $R$ with some
$(\eps,d,\delta)$-super-regular pairs. If a graph $H$ with
$\Delta(H)\leq\Delta$ is embeddable into $R(N)$ then it is already
embeddable into $G$.
\end{theorem}

Our proof is almost identical to the proof in \cite{KSSz3}. Of
course one difference is that we have to compute explicitly the
dependence between the parameters. Furthermore, this is also a
slight strengthening of the original statement as there can be a
small number of exceptional vertices which may have smaller degrees
($\delta$ may be much smaller than $d$). We note that the recent
``arrangeable'' Blow-up Lemma \cite{BKTW} is also quantitative, but
first of all the bound on $\eps$ is somewhat weaker and second it
does not allow for the strengthening mentioned above. However, in a
recent application \cite{GS} we needed precisely this strengthening.
We believe that this quantitative version of the Blow-up Lemma will
find other applications as well.

In Section 2 we give the embedding algorithm. In Section 3 we show
that the algorithm is correct.

\section{The algorithm}
The main idea of the algorithm is the following.
We embed the vertices of $H$ one-by-one by following
a greedy algorithm, which works smoothly
until there is only a small proportion of $H$ left,
and then it may get stuck hopelessly.
To avoid that, we will set aside a positive proportion
of the vertices of $H$ as buffer vertices.
Most of these buffer vertices will be embedded
only at the very end by using a K\"onig-Hall argument.

\subsection{Preprocessing}
We will assume that $|V(H)|=|V(G)|=|\cup_iV_i|=n=rN$. We will assume
for simplicity, that the density of every super-regular pair in $G$
is exactly $d$. This is not a significant restriction, otherwise we
just have to put everywhere the actual density instead of $d$.

We will use the following parameters:
$$\eps\ll\eps'\ll\eps''\ll \eps'''\ll d'''\ll d''\ll d'\ll d,$$
where $a\ll b$ means that $a$ is small enough compared to $b$. For
example we can select the parameters in the following explicit way:
$$d' = \frac{\delta d^{\Delta}}{8r\Delta},
d'' = (d')^3 = \left( \frac{\delta d^{\Delta}}{8r\Delta}
\right)^{3}, d''' = (d'')^2 = \left( \frac{\delta
d^{\Delta}}{8r\Delta} \right)^{6}, \eps''' = (d''')^2 = \left(
\frac{\delta d^{\Delta}}{8r\Delta} \right)^{12},$$
$$\eps'' = \left( \frac{\delta
d^{\Delta}}{8r\Delta} \right)^{2} (d''')^2 (\eps''')^3= \left(
\frac{\delta d^{\Delta}}{8r\Delta} \right)^{50}, \eps' = \left(
\frac{\delta d^{\Delta}}{8r\Delta} \right)^{2} (d''')^2 (\eps'')^3=
\left( \frac{\delta d^{\Delta}}{8r\Delta} \right)^{164} \;
\mbox{and}$$ $$ \eps = (\eps')^2 = \left( \frac{\delta
d^{\Delta}}{8r\Delta} \right)^{328}.$$

For easier reading, we will mostly use the letter $x$
for vertices of $H$,
and the letter $v$ for vertices of the host graph $G$.

Given an embedding of $H$ into $R(N)$, it defines an {\it assignment}
$$\psi:V(H)\to\{V_1,V_2,\ldots,V_r\},$$
and we want to find an {\it embedding}
$$\ffi:V(H)\to V(G),\quad\ffi\ \ \mbox{is one-to-one}$$
such that $\ffi(x)\in\psi(x)$ for all $x\in V(H)$.
We will write $X_i=\psi^{-1}(V_i)$ for $i=1,2,\dots,r$.

Before we start the algorithm, we order the vertices of $H$ into a
sequence $S=(x_1,x_2,\ldots,x_n)$ which is more or less, but not
exactly, the order in which the vertices will be embedded (certain
exceptional vertices will be brought forward). Let $m= r d'N$. For
each $i$, choose a set $B_i$ of $d'N$ vertices in $X_i$ such that
any two of these vertices are at a distance at least 3 in $H$. (This
is possible, for $H$ is a bounded degree graph.) These vertices
$b_1,\dots,b_m$ will be called the {\bf buffer vertices} and they
will be the last vertices in $S$.

The order $S$ starts with the neighborhoods
$N_H(b_1),N_H(b_2),\dots,N_H(b_m)$. The length of this initial
segment of $S$ will be denoted by $T_0$. Thus
$T_0=\ds\sum_{i=1}^m|N_H(b_i)|\leq \Delta m$.

The rest of $S$ is an arbitrary ordering
of the leftover vertices of $H$.

\subsection{Sketch of the algorithm}
In {\em Phase 1} of the algorithm we will embed the vertices in $S$
one-by-one into $G$ until all non-buffer vertices are embedded. For
each $x_j$ not embedded yet (including the buffer vertices) we keep
track of an ever shrinking host set $H_{t,x_j}$ that $x_j$ is
confined to at time $t$, and we only make a final choice for the
location of $x_j$ from $H_{t,x_j}$ at time $j$. At time 0,
$H_{0,x_j}$ is the cluster that $x_j$ is assigned to. For technical
reasons we will also maintain another similar set, $C_{t,x_j}$,
where we will ignore the possibility that some vertices are occupied
already. $Z_t$ will denote the set of occupied vertices. Finally we
will maintain a set $\mbox{Bad}_t$ of exceptional pairs of vertices.

In {\em Phase 2}, we embed the leftover vertices
by using a K\"onig-Hall type argument.

\subsection{Embedding Algorithm}

At time 0, set $C_{0,x}=H_{0,x}=\psi(x)$ for all $x\in V(H)$. Put
$T_1=d''n$.

\bs{\bf Phase 1.}\\
For $t\geq 1$, repeat the following steps.\\

{\em Step 1} (Extending the embedding).\quad We embed $x_t$.
Consider the vertices in $H_{t-1,x_t}$. We will pick one of these
vertices as the image $\ffi(x_t)$ by using the Selection Algorithm
(described below in Section \ref{selection}).\\

{\it Step 2} (Updating).\quad We set $$Z_t = Z_{t-1}\cup \{
\ffi(x_t)\},$$ and for each unembedded vertex $y$ (i.e. the set of
vertices $x_j,t<j\leq n$), set
$$C_{t,y}=\cases{C_{t-1,y}\cap N_G(\ffi(x_t))
&if $\{x_t,y\}\in E(H)$\cr
C_{t-1,y}&otherwise,}$$and
$$H_{t,y}=C_{t,y}\setminus Z_t.$$
We do not change the ordering at this step.\\

{\em Step 3} (Exceptional vertices in $G$).\quad

1. If $t\not\in \{1,T_0\}$, then go to Step 4.

2. If $t=1$, then we do the following (this is the part that is new
compared to the proof in \cite{KSSz3}). We find the 1st exceptional
set (denoted by $E_i^1$) consisting of those exceptional vertices
$v\in V_i$, $1\leq i\leq r$ for which there exists a $j\not= i$ such
that $(V_i,V_j)$ is $(\eps,d,\delta)$-super-regular, yet
$$deg_G(v,V_j) < (d-\eps) |V_j|.$$
(Note that $deg_G(v,V_j) \geq \delta |V_j|$ always holds by
super-regularity.) $\eps$-regularity implies that $|E_i^1|\leq r
\eps N$. We are going to change slightly the order of the vertices
in $S$. We choose a set $E_H^1$ of nonbuffer vertices $x\in H$ of
size $\sum_{i=1}^{r}|E_i^1|$ (more precisely $|E_i^1|$ vertices from
$X_i$ for all $1\leq i\leq r$) such that they are at a distance at
least 3 from each other. This is possible since $H$ is a bounded
degree graph and $\sum_{i=1}^r|E_i^1|$ is very small. We bring the
vertices in $E_H^1$ forward, followed by the remaining vertices in
the same relative order as before. For simplicity we keep the
notation $(x_1,x_2,\ldots,x_n)$ for the resulting order.
Furthermore, we slightly change the value of $T_0$ to
$T_0=\ds|E_H^1|+\sum_{i=1}^m|N_H(b_i)|$.\\

3. If $t=T_0$, then we do the following. We find the 2nd exceptional
set (denoted by $E_i^2$) consisting of those exceptional vertices
$v\in V_i$, $1\leq i\leq r$ for which $v$ is not covered yet in the
embedding and
$$\left|\left\{b:b\in B_i,v\in C_{t,b}\right\}\right|
<d''|B_i|.$$ Once again we are going to change slightly the order of
the remaining unembedded vertices in $S$. We choose a set $E_H^2$ of
unembedded nonbuffer vertices $x\in H$ of size
$\sum_{i=1}^{r}|E_i^2|$ (more precisely $|E_i^2|$ vertices from
$X_i$ for all $1\leq i\leq r$) with
$$H_{t,x}=H_{0,x}\setminus\{\ffi(x_j):j\leq t\}
=\psi(x)\setminus\{\ffi(x_j):j\leq t\}.$$ Thus in particular, if
$x\in X_i$, then $E_i^2\subset H_{t,x}$. Again we may choose the
vertices in $E_H$ as vertices in $H$ that are at a distance at least
3 from each other and any of the vertices embedded so far. We are
going to show later in the proof of correctness that this is
possible since $H$ is a bounded degree graph and
$\sum_{i=1}^r|E_i^2|$ is very small as well. We bring the vertices
in $E_H^2$ forward, followed by the remaining unembedded vertices in
the same relative order as before. Again, for simplicity we keep the
notation $(x_1,x_2,\ldots,x_n)$
for the resulting order.\\

{\em Step 4} (Exceptional vertices in $H$).\quad

1. If $T_1$ does not divide $t$, then go to Step 5.

2. If $T_1$ divides $t$, then we do the following. We find all
exceptional unembedded vertices $y\in H$ such that
$|H_{t,y}|\leq(d')^2n$. Once again we slightly change the order of
the remaining unembedded vertices in $S$. We bring these exceptional
vertices forward (even if they are buffer vertices), followed by the
non-exceptional vertices in the same relative order as before. Again
for simplicity we still use the notation $(x_1,x_2,\ldots,x_n)$ for
the new order. Note that it will follow from the proof, that if
$t\leq 2T_0$, then we do not find any exceptional vertices in $H$,
so we do not change the ordering at this step. \\

{\em Step 5} -
If there are no more unembedded non-buffer vertices left,
then set $T=t$ and go to Phase 2,
otherwise set $t\leftarrow t+1$ and go back to Step 1.

\bs{\bf Phase 2}\\
Find a system  of distinct representatives of the sets
$H_{T,y}$ for all unembedded $y$
(i.e. the set of vertices $x_j$, $T<j\leq n$).

\subsection{Selection Algorithm}\label{selection}

We distinguish two cases. Let $E_H = E_H^1 \cup E_H^2$.

{\bf Case 1.} $x_t \not\in E_H$.

We choose a vertex $v\in H_{t-1,x_t}$ as the image $\ffi(x_t)$ for
which the following hold for all unembedded $y$ with $\{x_t,y\}\in
E(H)$, \beq\label{1}(d-\eps)|H_{t-1,y}|\leq deg_G(v,H_{t-1,y})
\leq(d+\eps)|H_{t-1,y}|,\eeq \beq\label{2}(d-\eps)|C_{t-1,y}|\leq
deg_G(v,C_{t-1,y}) \leq(d+\eps)|C_{t-1,y}|\eeq and
\beq\label{3}(d-\eps)|C_{t-1,y}\cap C_{t-1,y'}| \leq
deg_G(v,C_{t-1,y}\cap C_{t-1,y'}) \leq(d+\eps)|C_{t-1,y}\cap
C_{t-1,y'}|,\eeq for at least a $(1-\eps')$ proportion of the
unembedded vertices $y'$ with $\psi(y')=\psi(y)$ and $\{y,y'\}
\not\in \mbox{Bad}_{t-1}$. Then we get $\mbox{Bad}_t$ by taking the
union of $\mbox{Bad}_{t-1}$ and the set of all of those pairs
$\{y,y'\}$ for which (\ref{3}) does not hold for $v=\ffi(x_t)$,
$C_{t-1,y}$ and $C_{t-1,y'}$. Thus note that we add at most $\Delta
\eps' N$ new pairs to $\mbox{Bad}_t$.

{\bf Case 2.} $x_t\in E_H$.

If $x_t\in X_i\cap E_H^l$, $l=1,2$, then we choose an arbitrary
vertex of $E_i^l$ as $\ffi(x_t)$. Note that for all $y\in N_H(x_t)$,
we have $C_{t-1,y}=\psi(y)$, \beq\label{4} deg_G(\ffi(x_t),
C_{t-1,y}) = deg_G(\ffi(x_t)) \geq \delta N = \delta
|C_{t-1,y}|,\eeq and \beq\label{5} deg_G(\ffi(x_t), H_{t-1,y}) \geq
deg_G(\ffi(x_t))-T_0-|E_H| \geq \delta N - 2 \Delta r d'N \geq
\frac{\delta}{2} N\eeq (using our choice of parameters). Here we
used super-regularity and the fact that $|E_H|\ll \Delta m$ which
will be shown later (Lemma \ref{E}).

\section{Proof of correctness}

The following claims state that our algorithm finds a good embedding
of $H$ into $G$.

\bigskip\noindent
{\bf Claim 1.}\ {\em Phase 1 always succeeds.}

\bigskip\noindent
{\bf Claim 2.}\ {\em Phase 2 always succeeds.}

\bigskip\noindent
If at time $t$, $S$ is a set of unembedded vertices $x\in H$ with
$\psi(x)=V_i$ (here and throughout the proof when we talk about time
$t$, we mean {\em after} Phase 1 is executed for time $t$, so for
example $x_t$ is considered embedded at time $t$), then we define
the bipartite graph $U_t$ as follows. One color class is $S$, the
other is $V_i$, and we have an edge between an $x\in S$ and a $v\in
V_i$ whenever $v\in C_{t,x}$.

In the proofs of the above claims the following lemma will play a
major role. First we prove the lemma for $t\leq T_0$, from this we
deduce that $|E_H|$ is small, then we prove the lemma for $T_0< t
\leq T$.

\begin{lemma}\label{egy}\pont
We are given integers $1\leq i\leq r$, $1\leq t\leq T_0$ and a set
$S\subset X_i$ of unembedded vertices at time $t$ with
$|S|\geq(d''')^2|X_i|=(d''')^2N$. If we assume that Phase 1
succeeded for all time $t'$ with $t'\leq t$, then apart from an
exceptional set $F$ of size at most $\eps''N$, for every vertex
$v\in V_i$ we have the following
$$deg_{U_t}(v)
=\left|\{x:x\in S,v\in C_{t,x}\}\right|\geq(1-\eps'')d(U_t)|S|
\quad\left(\geq\frac{d^{\Delta}}2|S|\right).$$
\end{lemma}

{\bf Proof.} In the proof of this lemma we will use the ``defect
form'' of the Cauchy-Schwarz inequality (just as in the original
proof of the Regularity Lemma \cite{SzemRegu}): if
$$\sum_{k=1}^m X_k=\frac{m}{n}\sum_{k=1}^n X_k
+D\qquad(m\leq n)$$ then
$$\sum_{k=1}^n X_k^2\geq\frac1n
\left(\sum_{k=1}^nX_k\right)^2+\frac{D^2n}{m(n-m)}.$$

Assume indirectly that the statement in Lemma \ref{egy} is not true,
that is, $|F|>\eps''N$. We take an $F_0\subset F$ with $|F_0|=\eps''
N$. Let us write $\nu(t,x)$ for the number of neighbors (in $H$) of
$x$ embedded by time $t$. Then in $U_t$ using the left side of
(\ref{2}) we get \beq\label{nagy}e(U_t)=d(U_t)|S||V_i| =\sum_{v\in
V_i}deg_{U_t}(v)=\sum_{x\in S}deg_{U_t}(x)$$ $$ =\sum_{x\in
S}|C_{t,x}|\geq\sum_{x\in S}(d-\eps)^{\nu(t,x)}N-\Delta r^2 \eps N^2
\geq (d-\eps)^{\Delta} |S|N - \Delta r^2 \eps N^2 \geq
\frac{d^{\Delta}}{2} |S|N,\eeq where the error term comes from the
neighbors of elements of $E_H^1$ (we are yet to start the embedding
of the vertices in $E_H^2$), since for them we cannot guarantee the
same lower bound.

We also have
$$\sum_{x\in S}\sum_{x'\in S}|N_{U_t}(x)\cap N_{U_t}(x')| =
\sum_{x\in S}\sum_{x'\in S}|C_{t,x}\cap C_{t,x'}|$$ $$ \leq
\sum_{x\in S}\sum_{x'\in S} (d+\eps)^{\nu(t,x)+\nu(t,x')}N + |S|N +
\Delta^2 |S| N + 2 \Delta r^2 \eps |S| N^2 + 2 \Delta \eps' N^3$$
\beq\label{7} \leq \sum_{x\in S}\sum_{x'\in S}
(d+\eps)^{\nu(t,x)+\nu(t,x')}N + 5 \Delta \eps' N^3.\eeq The error
terms come from the following $(x,x')$ pairs. For each such pair we
estimate $|C_{t,x}\cap C_{t,x'}|\leq N$. The first error term comes
from the pairs where $x=x'$. The second error term comes from those
pairs $(x,x')$ for which $N_H(x)\cap N_H(x')\not= \emptyset$. The
number of these pairs is at most $|S|\Delta(\Delta -1)\leq \Delta^2
|S|$. The third error term comes from those pairs $(x,x')$ for which
$x$ or $x'$ is a neighbor of an element of $E_H^1$. Finally we have
the pairs for which $\{x,x'\}\in \mbox{Bad}_t$. The number of these
pairs is at most $2t\Delta \eps' N \leq 2 \Delta \eps' N^2$.

Next we will use the Cauchy-Schwarz inequality with $m=\eps''N$ and
the variables $X_k, k=1, \ldots, N$ are going to correspond to
$deg_{U_t}(v)$, $v\in V_i$ (and the first $m$ variables to degrees
in $F_0$). Then we have
$$|D| = \eps'' \sum_{v\in V_i}deg_{U_t}(v) -  \sum_{v\in F_0}
deg_{U_t}(v)$$ \beq\label{D} \geq \eps'' \sum_{v\in V_i}deg_{U_t}(v)
- \eps''(1-\eps'') d(U_t)|S|N = (\eps'')^2 d(U_t)|S|N.\eeq

Then using (\ref{nagy}), (\ref{D}) and the Cauchy-Schwarz inequality
we get

$$\sum_{x\in S}\sum_{x'\in S}|N_{U_t}(x)\cap N_{U_t}(x')|
=\sum_{v\in V_i}(deg_{U_t}(v))^2$$
$$\geq\frac1N\left(\sum_{v\in V_i}deg_{U_t}(v)\right)^2
+(\eps'')^3d(U_t)^2N|S|^2$$
$$\geq\frac1N\left(\sum_{x\in S}
(d-\eps)^{\nu(t,x)}N - \Delta r^2 \eps N^2\right)^2
+(\eps'')^3d(U_t)^2N|S|^2$$
$$\geq\sum_{x\in S}\sum_{x'\in S}
(d-\eps)^{\nu(t,x)+\nu(t,x')}N - 2 \Delta \eps'
N^3+(\eps'')^3(d-\eps)^{2\Delta}N|S|^2,$$ which is a contradiction
with (\ref{7}), since $|S|\geq (d''')^2N$,
$$\left((d+\eps)^{\nu(t,x)+\nu(t,x')}-(d-\eps)^{\nu(t,x)+\nu(t,x')}\right)\ll
\Delta \eps,$$ and
$$ (\eps'')^3 (d-\eps)^{2\Delta} (d''')^2 \gg  \frac{d^{2\Delta}}{2} (d''')^2 (\eps'')^3  \geq \Delta \eps' \gg
\Delta \eps,$$ by the choice of the parameters. \qed

An easy consequence of Lemma \ref{egy} is the following lemma.

\begin{lemma}\label{E}\pont
In Step 3 we have $|E_i^2| \leq \eps '' N$ for every $1\leq i \leq
r$.
\end{lemma}

{\bf Proof.} Indeed applying Lemma \ref{egy} with $t=T_0$ and
$S=B_i$ (so we have $|S|=|B_i|= d' N > (d''')^2 N$) we get
$$(1-\eps'') d(U_t) |S| \geq \frac{d^{\Delta}}{2} |S| > d'' |S|,$$
and $E_i^2\subset F$. \qed

From this we can prove Lemma \ref{egy} for $t>T_0$ with $\eps'''$
instead of $\eps''$.

\begin{lemma}\label{egy'}\pont
We are given integers $1\leq i\leq r$, $T_0 < t\leq T$ and a set
$S\subset X_i$ of unembedded vertices at time $t$ with
$|S|\geq(d''')^2|X_i|=(d''')^2N$. If we assume that Phase 1
succeeded for all time $t'$ with $t'\leq t$, then apart from an
exceptional set $F$ of size at most $\eps'''N$, for every vertex
$v\in V_i$ we have the following
$$deg_{U_t}(v)
=\left|\{x:x\in S,v\in C_{t,x}\}\right|\geq(1-\eps''')d(U_t)|S|
\quad\left(\geq\frac{d^{\Delta}}2|S|\right).$$
\end{lemma}

{\bf Proof.} We only have to pay attention to the neighbors of the
elements of $E_H^2$, otherwise the proof is the same as the proof of
Lemma \ref{egy} with $\eps'''$ instead of $\eps''$. In (\ref{nagy})
the error term becomes $\Delta r \eps'' N^2$, coming from the
neighbors of elements of $E_H^2$.

In (\ref{7}) we have more bad pairs, namely all pairs $(x,x')$ where
$x$ or $x'$ is a neighbor of an element of $E_H^2$. These give an
additional error term of $2\Delta r \eps'' |S| N^2$. However, the
contradiction still holds, since
$$ (\eps''')^3 (d-\eps)^{2\Delta} (d''')^2 \gg  \frac{d^{2\Delta}}{2} (d''')^2 (\eps''')^3  \geq \Delta \eps'',$$
by the choice of the parameters. \qed

An easy consequence of Lemmas \ref{egy} and \ref{egy'} is the
following lemma.

\begin{lemma}\label{ketto}\pont
We are given integers $1\leq i\leq r$, $1\leq t\leq T$, a set
$S\subset X_i$ of unembedded vertices at time $t$ with $|S|\geq
d'''|X_i|=d'''N$ and a set $A\subset V_i$ with $|A|\geq
d'''|V_i|=d'''N$. If we assume that Phase 1 succeeded for all time
$t'$ with $t'\leq t$, then apart from an exceptional set $S'$ of
size at most $(d''')^2N$, for every vertex $x\in S$ we have the
following \beq\label{jolmetsz} |A\cap
C_{t,x}|\geq\frac{|A|}{2N}|C_{t,x}|. \eeq
\end{lemma}

{\bf Proof.} Assume indirectly that the statement is not true, i.e.
there exists a set $S'\subset S$ with $|S'|>(d''')^2N$ such that for
every $x\in S'$ (\ref{jolmetsz}) does not hold. Once again we
consider the bipartite graph $U_t=U_t(S',V_i)$. We have
$$\sum_{v\in A}deg_{U_t}(v)=\sum_{x\in S'}
|A\cap C_{t,x}|<\frac{|A|}{2N}\,\sum_{x\in S'}|C_{t,x}|
=\frac{|A|}{2N}\ d(U_t)|S'|N.$$ On the other hand, applying Lemmas
\ref{egy} or \ref{egy'} for $S'$ we get
$$\sum_{v\in A}deg_{U_t}(v)\geq(1-\eps'')d(U_t)|S'|(|A|-\eps''N)$$
contradicting the previous inequality. \qed

Finally we have
\begin{lemma}\label{harom}\pont
For every $1\leq t\leq T$ and for every vertex $y$ that is
unembedded at time $t$, if we assume that Phase 1 succeeded for all
time $t'$ with $t'\leq t$, then we have the following at time $t$
\beq\label{nagyH}|H_{t,y}|>d''N.\eeq
\end{lemma}

{\bf Proof.} We apply Lemma \ref{ketto} with $S_t$ the set of all
unembedded vertices in $X_i$ at time $t$, and $A_t=V_i\setminus Z_t$
(all uncovered vertices). Then for all but at most $(d''')^2N$
vertices $x\in S_t$ using (\ref{2}) and (\ref{4}) we get
\beq\label{11} |H_{t,x}| = |A_t \cap C_{t,x}| \geq \frac{|A_t|}{2N}
|C_{t,x}| \geq \frac{d'}{4} \delta (d-\eps)^{\Delta} N \geq (d')^2
N,\eeq if $|A_t|\geq (d'/2)N$. We will show next that in fact for
$1\leq t \leq T$, we have
$$|A_t|\geq |A_T|\geq (d'-d'') N \; \; \; \left( \geq \frac{d'}{2}
N\right),$$ so (\ref{11}) always holds. Assume indirectly that this
is not the case, i.e. there exists a $1\leq T'< T$ for which,
$$|A_{T'}| \geq (d'-d'')N \; \; \mbox{but} \; \; |A_{T'+1}| <
(d'-d'')N.$$ From the above at any given time $t$ for which $T_1|t$
and $1\leq t \leq T'$, in Step 4 we find at most $(d''')^2 N$
exceptional vertices in $X_i$. Hence, altogether we find at most
$$\frac{1}{d''} (d''')^2 N \ll d'' N$$
exceptional vertices in $X_i$ up to time $T'$. However, this implies
that at time $T'$ we still have many more than $(d'-d'')N$
unembedded buffer vertices in $X_i$, which in turn implies that
$|A_{T'+1}|\gg (d'-d'')N$, a contradiction. Thus we have
$$|A_T|\geq (d'-d'')N, \; \; T\leq rN - r d' N + r d'' N,$$
at time $T$ (or in Phase 2) we have at least $(d'-d'')N$ unembedded
buffer vertices in each $X_i$, and furthermore, for every $1\leq
t\leq T$ for all but at most $(d''')^2N$ vertices $x\in S_t$ we have
$$|H_{t,x}|> (d')^2 N.$$
Let us pick an arbitrary $1\leq t\leq T$ and an unembedded $y$ at
time $t$ (with $\psi(y)=V_i$). We have to show that (\ref{nagyH})
holds. Let $kd''n=kT_1 \leq t < (k+1)T_1$ for some $0\leq k \leq
T/T_1$. We distinguish two cases:

{\bf Case 1.} $y$ was not among the at most $(d''')^2N$ exceptional
vertices of $X_i$ found in Step 4 at time $kT_1$. Then
$$|H_{t,y}| \geq \left( \frac{\delta}{2} (d-\eps)^{\Delta} (d')^2 - r
d''\right) N.$$ Indeed, at time $kT_1$ we had $|H_{kT_1,y}|\geq
(d')^2N$. Until time $t$, $H_{t,y}$ could have been cut at most once
to a $\geq (\delta/2)$-fraction (if $y$ is a neighbor of an element
of $E_H$, there can be at most one such $E_H$-neighbor) and at most
$\Delta$ times to a $\geq (d-\eps)$-fraction (using (\ref{1}) and
(\ref{5})), and precisely $t-kT_1 \leq T_1 = r d'' N$ new vertices
were covered.

{\bf Case 2.} $y$ was among the at most $(d''')^2N$ exceptional
vertices of $X_i$ found in Step 4 at time $kT_1$. Then
$$|H_{t,y}| \geq \left( \frac{\delta}{2} (d-\eps)^{\Delta} (d')^2 - r
d'' - r (d''')^2 \right) N,$$ since at time $(k-1)T_1$ (we certainly
must have $k\geq 2$), $y$ was not exceptional, and because the
exceptional vertices were brought forward we have $t\leq kT_1 +
r(d''')^2N$. Thus in both cases we have $|H_{t,y}|> d'' N$, as
desired. \qed

Finally we show that the selection algorithm always succeeds in
selecting an image $\ffi(x_t)$.

\begin{lemma}\label{het}\pont For every $1\leq t \leq T$, if we
assume that Phase 1 succeeded for all time $t'$ with $t'\leq t$,
then Phase 1 succeeds for time $t$.
\end{lemma}

{\bf Proof.} We only have to consider Case 1 in the selection
algorithm. We choose a vertex $v\in H_{t-1,x_t}$ as the image
$\ffi(x_t)$ which satisfies (\ref{1}), (\ref{2}) and (\ref{3}). We
have by Lemma \ref{harom},
$$|H_{t-1,x_t}| \geq d'' N.$$
By $\eps$-regularity we have at most $2\eps N$ vertices in
$H_{t-1,x_t}$ which do not satisfy (\ref{1}) and similarly for
(\ref{2}). For (\ref{3}) we define an auxiliary  bipartite graph $B$
as follows. One color class $W_1$ is the vertices in $H_{t-1,x_t}$
and the other class $W_2$ is the sets $C_{t-1,y}\cap C_{t-1,y'}$ for
all pairs $\{y,y'\}$ where $\{x_t,y\}\in E(H)$, $\psi(y)=\psi(y')$,
and $\{y,y'\}\not\in \mbox{Bad}_{t-1}$. We put an edge between a
$v\in W_1$ and an $S\in W_2$ if inequality (\ref{3}) is not
satisfied for $v$ and $S$. Let us assume indirectly that we have
more than $\eps' N$ vertices $v\in W_1$ with $deg_B(v) > \eps'
|W_2|$. Then there must exist a $S\in W_2$ with
$$deg_B(S) > \eps' |W_1| \gg \eps N.$$
However, this is a contradiction with $\eps$-regularity since
$$|S| \geq (d-\eps)^{2\Delta} N \gg \eps N.$$
Here we used the fact that the pair corresponding to $S$ is not in
$\mbox{Bad}_{t-1}$. Thus altogether we have at most $4\eps N + \eps'
N \ll d'' N$ vertices in $H_{t-1,x_t}$ that we cannot choose and
thus the selection algorithm always succeeds in selecting an image
$\ffi(x_t)$, proving Claim 1. \qed

\bigskip\noindent
{\bf Proof of Claim 2.} We want to show that we can find a system of
distinct representatives of the sets $H_{T,x_j},T<j\leq n$, where
the sets $H_{T,x_j}$ belong to a given cluster $V_i$.

To simplify notation, let us denote by $Y$
the set of remaining vertices in $V_i$,
and by $X$ the set of remaining
unembedded (buffer) vertices assigned to $V_i$.
If $x=x_j\in X$ then write $H_x$ for its possible location
$H_{T,x_j}$ at time $T$. Also write $M=|X|=|Y|$.

The K\"onig-Hall condition for the existence of a system of distinct
representatives obviously follows from the following three
conditions: \beq\label{konig1} |H_x|>d'''M\quad\mbox{for all}\quad
x\in X, \eeq \beq\label{konig2} |\displaystyle\mathop\cup_{x\in
S}H_x|\geq(1-d''')M \text{for all subsets}S\subset X,|S|\geq d'''M,
\eeq \beq\label{konig3} |\displaystyle\mathop\cup_{x\in S}H_x|=M
\text{for all subsets}S\subset X,|S|\geq(1-d''')M. \eeq Equation
(\ref{konig1}) is an immediate consequence of Lemma \ref{harom},
(\ref{konig2}) is a consequence of Lemma \ref{egy}. Finally to prove
(\ref{konig3}), we have to show that every vertex in $Y\subset V_i$
belongs to at least $d''' |X|$ location sets $H_x$. However, this is
trivial from the construction of the embedding algorithm, in Step 3
of Phase 1 we took care of the small number of exceptional vertices
for which this is not true. This finishes the proof of Claim 2 and
the proof of correctness. \qed

\end{document}